\documentclass[12pt,twoside,openany,a4paper]{article}
\usepackage{amsmath}
\usepackage[usenames]{color}
\usepackage{mathrsfs}
\usepackage{amsfonts}
\usepackage{amssymb,amsmath}
\usepackage{cite}
\usepackage{cases}
\usepackage{amsthm}

\usepackage{hyperref}

\newtheorem{theo}{Theorem}[section]
\newtheorem{defi}[theo]{Definition}
\newtheorem{coro}[theo]{Corollary}
\newtheorem{lemm}[theo]{Lemma}
\newtheorem{propo}[theo]{Proposition}

\def\cl{\centerline}

\def\vs{\vspace*}
\def\ni{\noindent}
\def\Z{\mathbb{Z}}

\def\C{\mathbb{C}}

\begin{document}
  \cl{\textbf{{\large{Representations of Super $W(2,2)$ algebra $\mathfrak{L}$}}}}

 \cl{Hao Wang$^{1)}$, Huanxia Fa$^{2)}$, Junbo Li$^{2)}$}

\cl{\small $^{1)}$Wu Wen-Tsun Key Laboratory of Mathematics and School of Mathematical Sciences,}
\cl{\small University of Science and Technology of China, Hefei 230026, China}

\cl{\small $^{2)}$School of Mathematics and Statistics, Changshu Institute of Technology, Changshu 215500, China}

\vs{6pt}
\noindent{\bf{Abstract.}}
In paper, we study the representation theory of super $W(2,2)$ algebra ${\mathfrak{L}}$.
We prove that ${\mathfrak{L}}$ has no mixed irreducible modules and
give the classification of irreducible modules of intermediate series.
We determinate the conjugate-linear anti-involution of ${\mathfrak{L}}$ and give the unitary modules of intermediate series.


\noindent{\it{MR(2000) Subject Classification}: }\vs{18pt} 17B10, 17B65, 17B68.

\section{Introduction}
It is well known that the Virasoro algebra (named after the physicist  Miguel Angel Virasoro) is a very important infinite dimensional Lie algebra and is widely used in conformal field theory and string theory. After that much attention has been paid to the Virasoro type Lie algebras and superalgebras (which contains the Virasoro algebra as their subalgebras), including their constructions, structures and representations. The $W$-algebra $W(2,2)$ is certainly a Virasoro type Lie algebra, which plays important rolls in many areas of mathematics and physics (It was introduced by Zhang and Dong in \cite{ZD-CMP2008} for the study of classification of vertex operator algebras generated by vectors of weight $2$). It possesses a basis $\{L_{m},\,I_{m}|m\in\Z\}$ as a vector space over the complex field $\C$, with the Lie brackets $[L_{m},L_{n}]=(m-n)L_{m+n}$, $[L_{m},I_{n}]=(m-n)I_{m+n}$, $[I_{m},I_{n}]=0$. Structures and representations of $W(2,2)$ are extensively investigated in many references, such as \cite{CL-LA2012}, \cite{JP-JMP2010}, \cite{LGZ-JMP2008}, \cite{LS-AMSES2011}, \cite{LSX-CAM2008} and \cite{ZT-LMA2012}.

Some Lie superalgebras with $W$-algebra $W(2,2)$ as their even parts were constructed in \cite{WCB-JPAMT2012} as an application of the classification of Balinsky-Novikov super-algebras with dimension $2|2$.In this paper we consider the infinite dimensional Lie super  $W(2,2)$-algebra over the algebraic closed field  $\mathbb{C}$ (for convenience, we denote it  $\mathfrak{L}$), with the following non-vanishing brackets:
\begin{eqnarray}\label{brackets}
\begin{array}{lllll}
&[L_m,L_n]=(m-n)L_{m+n},&[L_m,I_n]=(m-n)I_{m+n},\vspace*{8pt}\\
&[L_m,H_n]=(\frac{m}{2}-n)H_{m+n},&[G_m,G_n]=I_{m+n},\vspace*{8pt}\\
&[L_m,G_n]=(\frac{m}{2}-n)G_{m+n},&[I_m,G_n]=(m-2n)H_{m+n}.
\end{array}
\end{eqnarray}
Finally we would like to make some remarks. We observe that many papers are forced on the Virasoro type Lie superalgebras which contain the super Virasoro algebra as their subalgebra (e.g., Refs. \cite{FLX-AC2011,YS-CS2009}), especially the $N=2$ super Virasoro algebras. It is easy to find that the algebra $\mathfrak{L}$ doesn't contain the super Virasoro Lie algebra as it's subalgebra, so the methods developed there are not applicable to $\mathfrak{L}$. All these make the study of $\mathfrak{L}$ more challengeable and attractive, we need to find some new methods to handel these problems.

\section{Weight module of ${\mathfrak{L}}$ with a finite dimensional weight space}
\subsection{Preliminaries and main results}
First we introduce some standard concepts.
Denote $\mathfrak{h}={\rm Span}_{\mathbb{C}}\{L_{0}\}$ the Cartan subalgebra of $\mathfrak{L}$.
The module of $\mathfrak{L}$ $M$ is called\ $\mathfrak{h}$-diagonalisable, if $M$ have the following decomposition:
\begin{equation}\label{WM-def-weight module}
M=\bigoplus_{\lambda\in\mathbb{C}}M_{\lambda},
\ \ \ M_{\lambda}=\{x\in M\,|\,L_{0}x=-\lambda\,x\},
\end{equation}
$M_{\lambda}$ is called the weight space of $\lambda$.
Denote ${\rm Supp}(M)=\{\lambda\in\mathbb{C}\,|\,M_{\lambda}\neq0\}$ the support set of $M$.
It is obvious£ºif $M$ is a irreducible weight module of $\mathfrak{L}$,
then there exists $\lambda\in\mathbb{C}$, such that ${\rm Supp}(M)\subset\lambda+\mathbb{Z}$.

The weight module $M$ of  $\mathfrak{L}$ is called  Harish-Chandra,
if all weight spaces of $M$ is finite dimensional.
The module  $M$ of  $\mathfrak{L}$ is called mixed, if there exists $\lambda,\,\mu\in{\rm Supp}(M)$,
such that ${\rm dim}M_{\lambda}=\infty$, ${\rm dim}M_{\mu}<\infty$.

Our main result of this section is the following theorem.

\begin{theo}\label{WM-theorem-main}
If $M$ is a irreducible weight module of $\mathfrak{L}$
and there exists $\lambda\in\mathbb{C}$, such that ${\rm dim}M_{\lambda}=\infty$, then ${\rm Supp}(M)=\lambda+\mathbb{Z}$,
and for each $k\in\mathbb{Z}$ we have ${\rm dim}M_{\lambda+k}=\infty$.
\end{theo}

From the above theorem, we immediately get the following two facts.

\begin{coro}
If $M$ is a irreducible weight module of $\mathfrak{L}$, with a finite dimensional weight space,
then $M$ ia a Harish-Chandra\ module.
\end{coro}

\begin{coro}
There exists no mixed module of $\mathfrak{L}$.
\end{coro}

\subsection{Proof of theorem \ref{WM-theorem-main}}
In this subsection,\ $M$\ is always a irreducible weight module of $\mathfrak{L}$. Notice that
\[\{L_{1},\,L_{2},\,I_{1},\,G_{1},\,H_{1},\,L_{-1},
\,L_{-2},\,I_{-1},\,G_{-1},\,H_{-1}\}\]
is a set of generators of $\mathfrak{L}$, we have the following proposition.

\begin{propo}\label{WM-principle fact}
If there exists $\mu\in\mathbb{C}$, $0\neq v\in M_{\mu}$ such that $L_{1}v=L_{2}v=I_{1}v=G_{1}v=H_{1}v=0$
or $L_{-1}v=L_{-2}v=I_{-1}v=G_{-1}v=H_{-1}v=0$, then $M$\ is a Harish-Chandra\ module.
\end{propo}

\begin{lemm}
If $\rm{dim} M_{\lambda}=\infty$, then there exists at most one $i\in\mathbb{Z}$,
such that $\rm{dim} M_{\lambda+i}<\infty$.
\end{lemm}

\ni{\it Proof.}
We prove it with reduction to absurdity. Without lose of generality, assume $\rm{dim} M_{\lambda+1}<\infty$,
$\rm{dim} M_{\lambda+j}<\infty$, in whick $j\in\mathbb{Z}$, $j\geq 3$. (If $j=2$, we immediately get a contradiction from \ref{WM-principle fact}).

Denote
\[V=(\cap_{A\in\{L, I, G, H\}}\ker\{A_{1}: M_{\lambda}\longrightarrow M_{\lambda+1}\})\cap (\cap_{A\in\{L, I, G, H\}}\ker\{A_{j}: M_{\lambda}\longrightarrow M_{\lambda+j}\}),\]
It is a subspace of $M_{\lambda}$.
$\rm{dim} M_{\lambda}=\infty$, $\rm{dim} M_{\lambda+1}<\infty$, together with
$\rm{dim} M_{\lambda+j}<\infty$ imply  $\rm{dim} V=\infty$.
 Act both sides of \eqref{brackets} of $V$, we get
\ \[L_{k}V=0\ \ \,\ \ \,k=1,\,j,\,j+1,\,j+2,\,\cdots\]
and
\ \[A_{k}V=0\ \ \,\ \ \,k=1,\,2,\,3,\,\cdots,\,\,A\in\{I,\,G,\,H\,\}.\]
If there exists $0\neq v\in V$ such that $L_{2}v=0$, Proposition \ref{WM-principle fact} implies $M$ is a Harish-Chandra module, contradicts with $\rm{dim} M_{\lambda}=\infty$\ .
Now we can assume for all  $0\neq v\in V$,  $L_{2}v\neq0$.
this implies $\rm{dim} L_{2}V=\infty$.

Denote $W=\cap_{A\in\{L,\,I,\,G,\,H\}}\ker\{A_{-1}:L_{2}V\longrightarrow M_{\lambda+1}\}$,
it is a subspace of $M_{\lambda+2}$. Since $\dim L_{2}V=\infty$, $\dim M_{\lambda+1}<\infty$,
we know $\rm{dim} W=\infty$,
this means $w=L_{2}v\ (0\neq v\in V)$ such that $A_{-1}w=0$,  in which $A\in\{\,L,\,I,\,G,\,H\,\}$.
For $k\geq j$\,$(\geq3)$, we have $L_{k}L_{2}v=(k-2)L_{k+2}v+L_{2}L_{k}v=0$.
So  $L_{-1}w=L_{k}w=0$ holds for all $k=1,\,j,\,j+1,\,j+2,\,j+3,\,\cdots$.
Since $[L_{-1},L_{k}]=-(1+k)L_{k-1}\neq 0$ for all $k\geq 1$,
it is obvious to see $L_{k}w=0$, $k=1,\,2,\,3,\,\cdots$,
This implies $A_{k}w=0$, $k=1,\,2,\,3,\,\cdots$, $A\in \{I,\,G,\,H\}$.
From Proposition \ref{WM-principle fact}, we get that $M$ is a Harish-Chandra module, contradicts with $\dim M_{\lambda}=\infty$.
\qed

Now we can assume for $i\in\mathbb{Z}^*$,  $\dim M_{\mu}<\infty$, $\dim M_{\mu+i}=\infty$.

\begin{lemm}\label{WM-lemma-002}
If $0\neq v\in M_{\mu-1}$ satisfies $L_{1}v=I_{1}v=G_{1}v=H_{1}v=0$,
then:

(1)\ \,$A_{k}v=0$.

(2)\ \,$A_{k}L_{2}v=0$, $\forall$ $k=1,\,2,\,3,\,\cdots,\,A\in\{I,\,G,\,H\}$.
\end{lemm}

\ni{\it Proof.}
(1) Since $[L_{1},A_{k}]=(\frac{1}{2}-k)A_{k+1}\neq0$ holds for $k=1,\,2,\,3,\,\cdots$, $A\in\{G,\,H\}$,
By induction on $k$, we have for $k=1,\,2,\,3,\,\cdots$,  $G_{k}v=H_{k}v=0$.
From $[G_{1},G_{1}]=I_{2}$ and $k=2,\,3,\,\cdots$, $[L_{1},I_{k}]=(1-k)I_{k+1}\neq0$,
Act both sides of the above equation on  $v$, we have $I_{k}v=0$ holds for $k=1,\,2,\,3,\,\cdots$.

(2)Holds immediately from (1) and $A_{k}L_{2}v=[A_{k},L_{2}]v+L_{2}A_{k}v$.
\qed

{\bf Proof of theorem \ref{WM-theorem-main}}\   Denote $V=\cap_{A\in\{L, I, G, H\}}\ker\{A_{1}: M_{\mu-1}\longrightarrow M_{\mu}\}$.
Since $\rm{dim} M_{\mu-1}=\infty$, $\rm{dim} M_{\mu}<\infty$,  $\rm{dim} V=\infty$.
For arbitrary $0\neq v\in V$, if $L_{2}v=0$, Proposition \ref{WM-principle fact} implies $M$ is a  Harish-Chandra module,
contradicts with the existence of a infinite weight space with $M$.

So for each $0\neq v\in V$, \ $L_{2}v\neq0$, this implies $\rm{dim} L_{2}V=\infty$.
Denote
\ \[W=\cap_{A\in\{L, I, G, H\}}\ker\{A_{-1}: L_{2}V\longrightarrow M_{\mu}\}.\]
Since $\rm{dim} L_{2}V=\infty$, $\rm{dim} M_{\mu}<\infty$, we have $\rm{dim} W=\infty$.
This implies there exists $w=L_{2}v$ ($0\neq v\in V$) such that for $A\in\{L,\,I,\,G,\,H\}$, $A_{-1}w=0$.
For $A\in\{G,\,H\}$ we have $[L_{-1},\,A_{k}]=(-\frac{1}{2}-k)A_{k-1}\neq 0$,
Act both sides of the above equation on $w$, we know for $k=-1,\,-2,\,-3,\,\cdots$,  $A_{k}w=0$.
Since for $k=-2,\,-3,\,-4,\,\cdots$, $[L_{-1},I_{k}]=(-1-k)I_{k-1}\neq 0$, and $[G_{-1},G_{-1}]=I_{-2}$,
Similar arguments as above, we have for $k=-1,\,-2,\,-3,\,\cdots$,  $I_{k}w=0$.

Act both sides of
\,\[[L_{-1},I_{1}]=-2I_{0},\ \ \ [L_{-1},G_{1}]=-\frac{3}{2}G_{0},\ \ \ [L_{-1},H_{1}]=-\frac{3}{2}H_{0},\]
on $w$, Lemma \ref{WM-lemma-002} implies for $k\in\mathbb{Z}, A\in\{I, G, H\}$,  $A_{k}w=0$.
Since $M$ is a irreducible weight module of $\mathfrak{L}$, we have\ $M=\mathfrak{U}(\mathfrak{L})w$
(In which $\mathfrak{U}(\mathfrak{L})$ is the universal enveloping algebra of $\mathfrak{L}$).
Denote $\textbf{L}$ the vector space generated by $\{L_{m}\,|\,m\in\mathbb{Z}\}$ over $\mathbb{C}$
(Obviously  $\textbf{L}$ is a centerless Virasoro Lie algebra ).
The above discussion implies $M=\mathfrak{U}(\textbf{L})w$,
which means $M$ is a irreducible weight module of centerless Virasoro algebra.
From the already known result: Virasoro owns no mixed weight weight module
(to be concrete, one can see \cite{MZ-JA2007}),
We have $M$ is a Harish-Chandra modules, contradicts with the assumption that $M$ has infinite dimensional weight space.
This complete the proof of our main theorem.
\qed

\section{Intermediate series module of ${\mathfrak{L}}$}

\subsection{Preliminaries and main results}
Denote $\mathfrak{h}={\rm Span}_{\mathbb{C}}\{L_{0}\}$ the Cartan subalgebra of $\mathfrak{L}$.
Weight module $M$ of $\mathfrak{L}$ is called $\mathfrak{h}$-diagonalisable , if $M$ has the decomposition:
\begin{equation}\label{def-weight module}
M=\bigoplus_{\lambda\in\mathbb{C}}M_{\lambda},
\ \ \ M_{\lambda}=\{x\in M\,|\,L_{0}x=-\lambda\,x\},
\end{equation}
$M_{\lambda}$ is a weight space of weight $\lambda$.
Denote ${\rm Supp}(M)=\{\lambda\in\mathbb{C}\,|\,M_{\lambda}\neq0\}$ the support set of $M$,
a $\mathfrak{L}$-module $M=M_{\overline{0}}\oplus M_{\overline{1}}$ is called the intermediate series module of $\mathfrak{L}$,
in which
\begin{equation}\label{def-weight module-evenoddpart}
M_{\overline{i}}=\bigoplus_{\lambda\in\mathbb{C}}M_{\overline{i}}^{\lambda},
\ \ \ M_{\overline{i}}^{\lambda}=\{x\in M\,|\,L_{0}x=-\lambda\,x\},
\end{equation}
if $M$ is irreducible and \rm{dim}$M_{\overline{i}}^{\lambda}\preceq 1$,
for $i=0$ or $1$.

First we recall some already known results about intermediate series modules of  Virasoro algebra. (To be concrete, one can see \cite{KS-MSRI1987}):
\begin{theo}\label{ISM-result for Vir}
 The intermediate series module $M$ of Virasoro algebra must be one the following modules\ $A_{a,b}$, $A(\alpha)$, $B(\beta)$,
  or their quotient modules:
  \begin{eqnarray*}
  A_{a,b}: &L_{i}u_{j}=(a-j+ib)u_{i+j},&\\
  A(\alpha): &L_{i}u_{j}=-(i+j)u_{i+j},\ j\neq0,& L_{i}u_{0}=-i(1+(i+1)\alpha)u_{i},\\
  B(\beta): &L_{i}u_{j}=-ju_{i+j},\ i+j\neq0,& L_{i}u_{-i}=i(1+(i+1)\beta)u_{0}.
\end{eqnarray*}
In which $i, j\in\mathbb{Z}$, $a, b, \alpha, \beta\in\mathbb{C}$.
$A_{a,b}$, $A(\alpha)$, $B(\beta)$\ as vector spaces over\ $\mathbb{C}$, own \{$u_{k}$ | $k\in\mathbb{Z}$\}\ as a basis.
\end{theo}

The main result of this section is the following theorem.
\begin{theo}\label{ISM-theorem main}
 If $M$ is a intermediate series module of $\mathfrak{L}$,
 it is also the  intermediate series module of Virasoro algebra.
  (In other words, $\forall k\in\mathbb{Z}$, $I_{k}, G_{k}, H_{k}$\ act trivally on $M$).
\end{theo}

\subsection{Proof of Theorem $\ref{ISM-theorem main}$}
We prove the theorem \ref{ISM-theorem main} by several lemmas. First we introduce some concepts.
Assume $M_{\overline{0}}$, $M_{\overline{1}}$ (modules of the  Virasoro algebra) are in forms of  $A_{a,b}$,
$M_{\overline{0}}=\rm{Span}_{\mathbb{C}}\{u_{k} \,|\, k\in\mathbb{Z}\}$,
$M_{\overline{1}}=\rm{Span}_{\mathbb{C}}\{v_{k} \,|\, k\in\mathbb{Z}\}$.
satisfy:
\begin{eqnarray}\label{ISM-eq-0}
  &L_{i}u_{j}=(a-j+ib)u_{i+j},& L_{i}v_{j}=(a'-j+ib')v_{i+j},\nonumber\\
  &I_{i}u_{j}=f(i,j)u_{i+j},& I_{i}v_{j}=\widetilde{f}(i,j)v_{i+j},\nonumber\\
  &G_{i}u_{i}=g(i,j)v_{i+j},& G_{i}v_{j}=\widetilde{g}(i,j)u_{i+j},\nonumber\\
  &H_{i}u_{j}=h(i,j)v_{i+j},& H_{i}v_{j}=\widetilde{h}(i,j)u_{i+j}.
\end{eqnarray}
in which $i, j\in\mathbb{Z}$, $a, a', b, b', f(i,j), \widetilde{f}(i,j), g(i,j),
\widetilde{g}(i,j), h(i,j), \widetilde{h}(i,j)\in\mathbb{C}$.

\begin{lemm}\label{ISM-Lemma I=0}
  $f(i,j)=\widetilde{f}(i,j)=0$, $\forall i, j\in\mathbb{Z}$.
\end{lemm}

\ni{\it Proof.}\ \ We prove this lemma case by case.

{\bf Case\ $1$}\ \ $a+mb\neq0$ $\forall m\in\mathbb{Z}$.

\ \ \ \ Act both sides of
\[
[L_{m},I_{n}]=(m-n)I_{m+n}.
\]
on\ $u_{t}$, we have:
\begin{eqnarray}\label{ISM-eq-1}
  &&(a-n-t+bm)f(n,t)v_{m+n+t}-(a-t+bm)f(n,m+t)v_{m+n+t}\nonumber\\
  &&=(m-n)f(m+n,t)v_{m+n+t},
\end{eqnarray}
This implies:
\begin{equation}\label{ISM-eq-2}
  (a-n-t+bm)f(n,t)-(a-t+bm)f(n,m+t)=(m-n)f(m+n,t).
\end{equation}
In \eqref{ISM-eq-2}, take $t=0$:
\begin{equation}\label{ISM-eq-3}
  (a+bm)f(n,m)=(a-n+bm)f(n,0)-(m-n)f(m+n,0).
\end{equation}
In \eqref{ISM-eq-2}, take $m=n$:
\begin{equation}\label{ISM-eq-4}
  (a+(b-1)m-t)f(m+t)-(a-t+bm)f(m,m+t)=0.
\end{equation}
In \eqref{ISM-eq-4}, take $t=0$:
\begin{equation}\label{ISM-eq-5}
  (a+bm)f(m,0)=(a+(b-1)m)f(m,0).
\end{equation}
In \eqref{ISM-eq-4}, take $t=-m$:
\begin{equation}\label{ISM-eq-6}
  (a+bm)f(m,-m)=(a+(b+1)m)f(m,0).
\end{equation}
In \eqref{ISM-eq-3}, take $m=-n$:
\begin{equation}\label{ISM-eq-7}
  (a-bn)f(n,-n)=(a-(b+1)n)f(n,0)+2nf(0,0).
\end{equation}
Associate \eqref{ISM-eq-6} with \eqref{ISM-eq-7}:
\begin{equation}\label{ISM-eq-8}
  f(m,0)=\frac{a+bm}{a}f(0,0),\ m\neq0.
\end{equation}
In \eqref{ISM-eq-2}, take $t=m, n=0$:
\begin{equation}\label{ISM-eq-9}
  (a-m+bm)(f(0,m)-f(0,2m))=mf(m,m).
\end{equation}
In \eqref{ISM-eq-3}, take $n=0$:
\begin{equation}\label{ISM-eq-10}
  (a+bm)f(0,m)=(a+bm)f(0,0)-mf(m,0).
\end{equation}
Associate \eqref{ISM-eq-8} with \eqref{ISM-eq-10}:
\begin{equation}\label{ISM-eq-11}
  f(0,m)=\frac{a-m}{a}f(0,0).
\end{equation}
From \eqref{ISM-eq-9} and \eqref{ISM-eq-11}:
\begin{equation}\label{ISM-eq-12}
  f(m,m)=\frac{a-m+bm}{a}f(0,0).
\end{equation}
Take \eqref{ISM-eq-8} together with \eqref{ISM-eq-11},  \eqref{ISM-eq-3} becomes:
\begin{equation}\label{ISM-eq-13}
  f(n,m)=\frac{f(0,0)(a+bn-m)}{a}.
\end{equation}
holds for all $m,n\in\mathbb{Z}$.
Act both sides of $[I_{m},I_{n}]=0$ on $u_{k}$,
we have $\frac{f(0,0)}{a}=0$.
From the above discussion we know for all $m,n\in\mathbb{Z}$, $f(m,n)=0$.

{\bf Case\ $2$} $a\neq0$, $0\in\{a+bn \,|\, n\in\mathbb{Z}\}$.

Under this assumption we know $b\neq0, -1$ and there exists  $p\in\mathbb{Z}^{*}$ such that $a=bp$.
Similar arguments as in case $1$: $f(m,n)=0$ $\forall n\neq-p$.
together with\ $\eqref{ISM-eq-4}$:
\begin{equation}\label{ISM-eq-14}
 f(0,-p)=0
\end{equation}
In \eqref{ISM-eq-2} take $t=-p, m=0$, together with \eqref{ISM-eq-14}:
$f(n,-p)=0$ $\forall n\neq0$.
From the above discussion, we know: $f(m,n)=0$ $\forall m, n\in\mathbb{Z}$.

{\bf Case $3$} $a\in\mathbb{Z}$.

Since $A_{a,b}\cong A_{0,b}$, without lose of generality, we assume $a=0$.

{\bf Case $3.1$} $b\neq0, -1$. Similar arguments as in $1$, we have:
$f(m,n)=0$ holds $\forall m, n\in\mathbb{Z}$.

{\bf Case $3.2$} $b=-1$.

In this situation, $\eqref{ISM-eq-2}$ and $\eqref{ISM-eq-3}$ turn into the following formula:
\begin{equation}\label{ISM-eq-15}
  (m+n+t)f(n,t)-(m+t)f(n,m+t)=(n-m)f(m+n,t).
\end{equation}
\begin{equation}\label{ISM-eq-16}
  mf(n,m)=(m+n)f(n,0)+(m-n)f(m+n,0).
\end{equation}
In \eqref{ISM-eq-15} take $n=t=0$:
\[
mf(0,0)-mf(0,m)=-mf(m,0).
\]
this implies:
\begin{equation}\label{ISM-eq-17}
  f(0,m)=f(m,0)+f(0,0),\ m\neq0.
\end{equation}
In \eqref{ISM-eq-15} take $t=m, m=-n$:
\begin{equation}\label{ISM-eq-18}
  mf(n,m)=(m-n)f(n,m-n)+2nf(0,m).
\end{equation}
In \eqref{ISM-eq-16}\ replace $m$ with $m-n$:
\begin{equation}\label{ISM-eq-19}
  (m-n)f(n,m-n)=mf(n,0)+(m-2n)f(m,0).
\end{equation}
Take together \eqref{ISM-eq-16}--\eqref{ISM-eq-19}, we have:
\begin{equation}\label{ISM-eq-20}
  (m-n)f(m+n,0)=mf(m,0)-nf(n,0)+2nf(0,0).
\end{equation}
In \eqref{ISM-eq-20} take $m=n\neq0$, we get $f(0,0)=0$.
$\eqref{ISM-eq-20}$ becomes:
\begin{equation}\label{ISM-eq-21}
  (m-n)f(m+n,0)=mf(m,0)-nf(n,0).
\end{equation}
By induction:
\begin{equation}\label{ISM-eq-22}
  f(m,0)=c+dm,\ \mbox{ÆäÖÐ}\ c, d\in\mathbb{C}.
\end{equation}
Associative \eqref{ISM-eq-20} with \eqref{ISM-eq-16}:
\begin{equation}\label{ISM-eq-23}
  f(n,m)=f(n,0)+f(m,0),\ m\neq0.
\end{equation}
From \eqref{ISM-eq-22} and \eqref{ISM-eq-23}:
\begin{equation}\label{ISM-eq-24}
  f(n,m)=2c+d(m+n), \ m\neq0.
\end{equation}
Replaced in \eqref{ISM-eq-15} we have $c=0$, in other words: $f(m,n)=d(m+n)$.
Together with $[I_{m},I_{n}]=0$ we know:
 $f(m,n)=0$ $\forall m, n\in\mathbb{Z}$.

{\bf Case $3.3$} $b=0$.

In this situation \eqref{ISM-eq-2} becomes:
\begin{equation}\label{ISM-eq-25}
  (-n-t)f(n,t)+tf(n,m+t)=(m-n)f(m+n,t).
\end{equation}
In \eqref{ISM-eq-25}, take $t=0, m=-n$:
\begin{equation}\label{ISM-eq-26}
  f(n,0)=2f(0,0),\ n\neq0.
\end{equation}
In \eqref{ISM-eq-25}, take $t=0$:
\begin{equation}\label{ISM-eq-27}
  f(k,0)=0,\ \forall\, k\in\mathbb{Z}.
\end{equation}
In \eqref{ISM-eq-25}, take $t=-m$, together with \eqref{ISM-eq-27}:
\begin{equation}\label{ISM-eq-28}
 f(n,-m)=f(m+n,-m),\ m\neq n.
\end{equation}
In \eqref{ISM-eq-25}, take $t=1$:
\begin{equation}\label{ISM-eq-29}
  (-n-1)f(n,1)+f(n,m+1)=(m-n)f(m+n,1).
\end{equation}
In \eqref{ISM-eq-29}, take $m=-1$, associate with \eqref{ISM-eq-26} we get:
\begin{equation}\label{ISM-eq-30}
  f(n,1)=f(n-1,1),\ n\neq-1.
\end{equation}
By induction on \eqref{ISM-eq-30}, we have:
\begin{equation}\label{ISM-eq-31}
  f(n,1)=
  \begin{cases}
    c_{1}, &n\geq-1,\\
    c_{2}, &n\leq-2.
  \end{cases}
\end{equation}
 \eqref{ISM-eq-28} and \eqref{ISM-eq-31} imply $c_{1}=c_{2}$.
From \eqref{ISM-eq-30} we have $f(n,m)=mc$ holds for all $m, n\in\mathbb{Z}$.
Act both sides of $[I_{m,I_{n}}]=0$ on $u_{k}$, we get $c=0$,
This implies for each $m, n\in\mathbb{Z}$,  $f(m,n)=0$.

Similarly, we have:
 $\widetilde{f}(m,n)=0$ $\forall m, n\in\mathbb{Z}$.
\qed

Lemma \ref{ISM-Lemma I=0}\ and the Lie bracket $[I_{m},G_{n}]=H_{m+n}$ imply the following lemma.

\begin{lemm}\label{ISM-Lemma H=0}
   $h(m,n)=\widetilde{h}(m,n)=0$,
   $\forall m, n\in\mathbb{Z}$.
\end{lemm}

Lemma \ref{ISM-Lemma I=0}\ and the Lie bracket $[G_{m},G_{n}]=I_{m+n}$ imply the following lemma.
\begin{lemm}\label{ISM-Lemma G-pre}
  \[
  g(m,k)\widetilde{g}(n,m+k)+g(n,k)\widetilde{g}(m,n+k)=0.
  \]
  $\forall m, n, k\in\mathbb{Z}$

\end{lemm}

\begin{lemm}\label{ISM-Lemma G-Main}
 At least one of $g(m,n)$ $\widetilde{g}(m,n)$ is $0$.
\end{lemm}

\ni{\it Proof.} Suppose not, assume neither  $g(m,n)$ nor $\widetilde{g}(m,n)$ is $0$.

Act both sides of $[L_{0},G_{n}]=-nG_{n}$ on $u_{k}$, we know:
\[
 (a'-a)g(n,k)=0,
\]
 $\forall n, k\in\mathbb{Z}$.
Under the assumption  $g(m,n)\not\equiv0$, we get $a=a'$.
Act both sides of $[L_{i},G_{j}]=(\frac{i}{2}-j)G_{i+j}$ on $u_{k}$:
\begin{equation}\label{ISM-eq-32}
 (\frac{i}{2}-j)g(i+j,k)=(a-(k+j)+ib')g(j,k)-(a-k+ib)g(j,i+k).
\end{equation}
In \eqref{ISM-eq-32}, take $i=2j$:
\begin{equation}\label{ISM-eq-33}
  (a-k+2jb)g(j,k+2j)=(a-(k+j)+2jb')g(j,k).
\end{equation}
In \eqref{ISM-eq-32} take $i=-2j$:
\begin{equation}\label{ISM-eq-34}
  -2jg(-j,k)=(a-(k+j)-2jb')g(j,k)-(a-k-2jb)g(j,k-2j).
\end{equation}
\ \eqref{ISM-eq-34} times $a-(k-j)+2jb'$,
and replace the last item of \eqref{ISM-eq-34} by \eqref{ISM-eq-33}:
\begin{eqnarray}\label{ISM-eq-35}
  &&-2j(a-(k-j)+2jb')g(-j,k)\nonumber\\
  &&\ \ =(a-(k-j)+2jb')(a-(k+j)-2jb')g(j,k)\nonumber\\
  &&\ \ \ \ \ \ -(a-k-2jb)(a-(k-2j)+2jb)g(j,k)\nonumber\\
  &&\ \ =2j(a-(k-j)+2jb'+2jt)g(j,k).
\end{eqnarray}
In which $t=b'^{2}-{(b+\frac{1}{2})}^{2}$.
Similarly: for \eqref{ISM-eq-32}, let $j=-j, i=2j$, $j=-j, i=-2j$:
\begin{equation}\label{ISM-eq-36}
  2j(a-(k+j)-2jb')g(j,k)=-2j(a-(k+j)-2jb'-2jt)g(-j,k).
\end{equation}
Together with \eqref{ISM-eq-35} and \eqref{ISM-eq-36}:
\begin{eqnarray}\label{ISM-eq-37}
  &\ &((a-(a+j)-2jb')(a-(k-j)+2jb')-((a-k-j-2jb')-2jt)\nonumber\\
  &\ \ \ \ \ &   \times((a-k+j+2jb')+2jt))g(j,k)=0.
\end{eqnarray}
From \eqref{ISM-eq-37} we have:
\begin{equation}\label{ISM-eq-38}
  4j^{2}t(t+2b'+1)g(j,k)=0.
\end{equation}
 \eqref{ISM-eq-32} implies: there sxists\ $k_{0}\in\mathbb{Z}$
such that $g(1,k_{0})\neq0$
(If not,  \eqref{ISM-eq-32} implies $g(m,n)\equiv0$, contradicts with our assumption).
Together with \eqref{ISM-eq-38}, we have:
\begin{equation}\label{ISM-eq-39}
      \begin{cases}
        b'=b+\frac{1}{2},\ \mbox{or}\\
        b'=-(b+\frac{1}{2}),\ \mbox{or}\\
        b'=-b-\frac{3}{2},\ \mbox{or}\\
        b'=b-\frac{1}{2}.
      \end{cases}
\end{equation}
In the following, we determine $g(m,n), \widetilde{g}(m,n)$ case by case.

{\bf Case\ $1$}\ \ $b'=b+\frac{1}{2}$.

Assume $a-k+2b\neq0$, $\forall k\in\mathbb{Z}$.
In \eqref{ISM-eq-33}, take $j=1$:
\begin{equation}\label{ISM-eq-40}
  g(1,k)=
  \begin{cases}
    x_{0},& k\ \mbox{is even},\\
    x_{1},& k\ \mbox{is odd}.
  \end{cases}
\end{equation}
In \eqref{ISM-eq-32}, take $j=1$,  $i, k$ are even integers:
\begin{eqnarray}\label{ISM-eq-41}
  (\frac{i}{2}-1)g(i+1,k)\!\!\!&=&\!\!\!(a-(k+1)+ib')g(1,k)-(a-k+ib)g(1,k+i)\nonumber\\
                         \!\!\!&=&\!\!\!(a-(k+1)+ib')x_{1}-(a-k+ib)x_{0}.
\end{eqnarray}
Associative \eqref{ISM-eq-33} with $b'=b+\frac{1}{2}$:
\begin{eqnarray}\label{ISM-eq-42}
  (\frac{i}{2}-1)g(i+1,k)\!\!\!&=&\!\!\!(\frac{i}{2}-1)g(i+1,k+2(i+1))\nonumber\\
                         \!\!\!&=&\!\!\!(a-(k+(2i+1)+1)+ib')x_{1}\nonumber\\
                         &\ \ &-(a-(k+(2i+1))+ib)x_{0}.
\end{eqnarray}
Obviously we have $x_{0}=x_{1}$.
Similar arguments following \eqref{ISM-eq-38}, we have:
\begin{equation}\label{ISM-eq-43}
  g(i,j)=d_{1},\ \forall\, i, j\in\mathbb{Z}.
\end{equation}
In which $d_{1}\in\mathbb{C}^{*}$ is a constant.

Now we assume  $a-k'+2b=0$, $\forall k'\in\mathbb{Z}$.
 Together with \eqref{ISM-eq-33} we know: $(a-k+2b)g(1,k)=(a-k+2b)g(1,k+2)$.
This implies:
\begin{equation}\label{ISM-eq-44}
  g(1,k)=
  \begin{cases}
    x_{0},\ k>k',\ k\ \mbox{is even},\\
    x_{1},\ k>k',\ k\ \mbox{is odd}.
  \end{cases}
\end{equation}

\begin{equation}\label{ISM-eq-45}
  g(1,k)=
  \begin{cases}
    y_{0},\ k\preceq k',\ k\ \mbox{is even},\\
    y_{1},\ k\preceq k',\ k\ \mbox{is odd}.
  \end{cases}
\end{equation}
Together with \eqref{ISM-eq-32}:
\begin{equation}\label{ISM-eq-46}
  g(i,k)=
  \begin{cases}
    x_{0}, k>k',\ k+i-1>k',\ k,\ i-1\ \mbox{is even},\\
    x_{1}, k>k',\ k+i-1>k',\ k,\ i-1\ \mbox{is odd}.
  \end{cases}
\end{equation}

\begin{equation}\label{ISM-eq-47}
  g(i,k)=
  \begin{cases}
    y_{0}, k\preceq k',\ k+i-1\preceq k',\ k,\ i-1\ \mbox{is even},\\
    y_{1}, k\preceq k',\ k+i-1\preceq k',\ k,\ i-1\ \mbox{is odd}.
  \end{cases}
\end{equation}
Select some $k, j\in\mathbb{Z}$, satisfy
$a-k+2jb\neq0, k\preceq k', k+j-1\preceq k', k+2j>k', k+3j-1>k'$,
one of  $k, j$ is even and the other is odd.
From \eqref{ISM-eq-33} we know:
\[
 (a-k+2jb)g(j,k)=(a-k+2jb)g(j,k+2j).
\]
this implies: $x_{0}=y_{0}, x_{1}=y_{1}$.
Similar arguments, we have:
$g(m,n)=d_{1}$, $\forall m, n\in\mathbb{Z}$, .
In which $d_{1}\in\mathbb{C}^{*}$.

{\bf Case $2$}\ \ $b'=-(b+\frac{1}{2})$.

In \eqref{ISM-eq-32}, take $i=2j$:
\begin{equation}\label{ISM-eq-48}
  (a-k-2j-2jb)g(j,k)=(a-k+2jb)g(j,k+2j).
\end{equation}
In \eqref{ISM-eq-32}, take $i=-2j$:
\begin{equation}\label{ISM-eq-49}
  -2jg(-j,k)=(a-k+2jb)g(j,k)-(a-k-2jb)g(j,k-2j)=-2jg(j,k).
\end{equation}
The second equation on the right side of \eqref{ISM-eq-49} follows from \eqref{ISM-eq-48} with the replacement $kj$ is replaced by  $k-2$.
This implies: $g(j,k)=g(-j,k)$ holds $\forall\, j, k\in\mathbb{Z}$.
Together with \eqref{ISM-eq-32}:
\begin{eqnarray}\label{ISM-eq-50}
  &&(a-(k+j)+ib')g(j,k)-(a-k+ib)g(j,k+i)\nonumber\\
  &&=(\frac{i}{2}-j)g(i+j,k)\nonumber\\
  &&=(\frac{i}{2}-j)g(-i-j,k)\nonumber\\
  &&=-(a-(k-j)-ib')g(-j,k)+(a-k-ib)g(-j,k-i).
\end{eqnarray}
This implies:
\begin{equation}\label{ISM-eq-51}
  (a-k-ib)g(j,k-i)-2(a-k)g(j,k)+(a-k+ib)g(j,k+i)=0.
\end{equation}
In \eqref{ISM-eq-48}, take $j=1$ and replace $k$ by $k+2$:
\begin{equation}\label{ISM-eq-52}
  g(1,k+4)=\frac{(a-k-4-2b)(a-k-2-2b)}{(a-k-2+2b)(a-k+2b)}g(1,k).
\end{equation}
Similarly, we get an expression of  $g(1,k-4)$ in form of \eqref{ISM-eq-52}.
IN \eqref{ISM-eq-51}, take $i=4$, together with \eqref{ISM-eq-52}:
\begin{eqnarray}\label{ISM-eq-53}
  &&(\frac{(a-k+4+2b)(a-k+2+2b)}{(a-k+2-2b)(a-k-2b)}\nonumber\\
  &&\times(a-k+4b)-2(a-k)+(a-k+4b)\nonumber\\
  &&\times\frac{(a-k-4-2b)(a-k-2-2b)}{(a-k-2+2b)(a-k+2b)})g(1,k)=0.
\end{eqnarray}
 From \eqref{ISM-eq-48} and the similar arguments following \eqref{ISM-eq-38},
we know the coefficient of $g(1,k)$ is $0$.
This implies:  either $b=-1$ or $-\frac{1}{2}$.
Since $b=-\frac{1}{2}, b'=0$ is just in case $1$,
we only need to consider $b=-1, b'=\frac{1}{2}$.
In \eqref{ISM-eq-33}, take $j=1$,
we know $(a-k)g(1,k)$ $k$ is a constant $\forall k$.
Assume:
\begin{equation}\label{ISM-eq-54}
  (a-k)g(1,k)=
  \begin{cases}
    x_{0},\ k\ \mbox{is even},\\
    x_{1},\ k\ \mbox{is odd}.
  \end{cases}
\end{equation}
In \eqref{ISM-eq-51}, take $i=j=1$, we have $x_{0}=x_{1}$.
In other words, $(a-k)g(1,k)$ is a constant.
If for $k_{1}\in\mathbb{Z}$, $a-k_{1}=0$ holds,
then $\forall k\neq k_{1}$:
\[
 g(1,k)=\frac{(a-k_{1})}{a-k}g(1,k_{1})=0.
\]
This ensures $g(m,n)\equiv0$, contradicts with our assumption.
From the discussion above, we have $(a-k)g(1,k)=d_{2}$, in which $d_{2}\in\mathbb{C}^{*}$.
In \eqref{ISM-eq-32}, take $j=1$:
\begin{equation}\label{ISM-eq-55}
  (\frac{i}{2}-1)g(i+1,k)=(a-(k+1)+\frac{i}{2})g(1,k)-(a-k-1)g(1,k+1).
\end{equation}
this implies $g(i,k)=g(1,k)$ holds for all $k\in\mathbb{Z}, k\neq3$.
In \eqref{ISM-eq-32}, take $j=2, i=1$:
\begin{equation}\label{ISM-eq-56}
  -\frac{3}{2}g(3,k)=(a-(k+2)+\frac{1}{2})g(2,k)-(a-k-1)g(2,k+1).
\end{equation}
this implies: $g(3,k)=g(1,k)$.
Now, we have $g(i,j)=\frac{d_{2}}{a-j}$, in which $d_{2}\in\mathbb{C}$.

{\bf Case\ $3$}\ \ $b'=-b-\frac{3}{2}$.

Similar arguments as in case $2$, we have $2jg(j,k)=2jg(j,k-2j)$, together with:
\begin{equation}\label{ISM-eq-57}
  (a-k-j+ib')g(j,k)+(a-k-j-2i-ib')g(j,k+2j)=2(a-k-i-j)g(j,k+i).
\end{equation}
\begin{eqnarray}\label{ISM-eq-58}
  &&(\frac{(a-k-3-2b')(a-k-5-2b')}{(a-k-1+2b')(a-k-3+2b')}\nonumber\\
  &&\times(a-k-1+4b')-2(a-k-5)+(a-k-9-4b')\nonumber\\
  &&\times\frac{(a-k-7+2b')(a-k-5+2b')}{(a-k-9-2b')(a-k-7-2b')})g(1,k+4)=0.
\end{eqnarray}
This implies $b=-\frac{3}{2}$ or $-1$.
The case $b=-1, b'=-\frac{1}{2}$ is just in case $1$.
So we only need to consider $b=-\frac{3}{2}, b'=0$.
In \eqref{ISM-eq-33}, take $j=1$:
\begin{equation}\label{ISM-eq-59}
  (a-k-1)g(1,k)=(a-k-3)g(1,k+2).
\end{equation}
Similar arguments as in case $2$, we know: $(a-k-1)g(1,k)$ is a constant.
If there exists $k_{1}\in\mathbb{Z}$ such that $a-K_{1}-1=0$,
then for all $k\neq k_{1}$, we have:
\[
 g(1,k)=\frac{a-k_{1}-1}{a-k-1}g(1,k_{1})=0.
\]
this ensures $g(m,n)\equiv0$, contradicts with our assumption.
Therefor, $g(1,k)=\frac{d_{3}}{a-k-1}$, in which $d_{3}\in\mathbb{C}^{*}$.
In \eqref{ISM-eq-32}, take $j=1$:
\begin{equation}\label{ISM-eq-60}
  (\frac{i}{2}-1)g(i+1,k)=(a-k-1)g(1,k)-(a-k-i\frac{3}{2})g(1,k+i).
\end{equation}
This implies:
\begin{equation}\label{ISM-eq-61}
  g(i,k)=g(1,k+i-1)=\frac{d_{3}}{a-k-i},\ i\neq3.
\end{equation}
In \eqref{ISM-eq-32}, take $i=1, j=2$:
\begin{equation}\label{ISM-eq-62}
  -\frac{3}{2}g(3,k)=(a-k-2)g(2,k)-(a-k-\frac{3}{2})g(2,k+1)=-\frac{3}{2}g(2,k+1).
\end{equation}
This implies:
\begin{equation}\label{ISM-eq-63}
  g(i,j)=\frac{d_{3}}{a-i-j}
\end{equation}
in which $d_{3}\in\mathbb{C}^{*}$.

{\bf Case\ $4$}\ \ $b'=b-\frac{1}{2}$.

Similar as in case $1$, we have: $\widetilde{g}(m,n)=d'_{1}$.

Take together all the cases $1--4$ and make the similar argument for $\widetilde{g}(m,n)$, we have:
\begin{equation}\label{ISM-eq-64}
 g(i,j)=
  \begin{cases}
    d_{1}, &b'=b+\frac{1}{2},\\
    \frac{d_{2}}{a-j}, &b'=\frac{1}{2},\ b=-1,\\
    \frac{d_{3}}{a-i-j}, &b'=0,\ b=-\frac{3}{2},\\
    g(i,j), &b'=b-\frac{1}{2}.
  \end{cases}
\end{equation}
\begin{equation}\label{ISM-eq-65}
  \widetilde{g}(i,j)=
  \begin{cases}
    d'_{1}, &b=b'+\frac{1}{2},\\
    \frac{d'_{2}}{a-j}, &b=\frac{1}{2},\ b'=-1,\\
    \frac{d'_{3}}{a-i-j}, &b=0,\ b'=-\frac{3}{2},\\
    \widetilde{g}(i,j), &b=b'-\frac{1}{2}.
  \end{cases}
\end{equation}
In which\ $i, j\in\mathbb{Z}, d_{1}, d'_{1}, d_{2}, d'_{2}, d_{3}, d'_{3}\in\mathbb{C}^{*}$.
Obviously only the two cases are possible:
\begin{equation}\label{ISM-eq-66}
  \begin{cases}
    g(i,j)=d_{1},\ \widetilde{g}(i,j)=\widetilde{g}(i,j), &b'=b+\frac{1}{2},\\
    g(i,j)=g(i,j),\ \widetilde{g}(i,j)=d'_{1}, &b'=b-\frac{1}{2}.
  \end{cases}
\end{equation}
associate with \eqref{ISM-eq-66}\ and lemma\ \ref{ISM-Lemma G-pre}, we have:
either $g(m,n)\equiv0$ or $\widetilde{g}(m,n)\equiv0$, this contradicts with our assumption.
This complete the proof of our lemma.
\qed

For $M_{\overline{0}}, M_{\overline{1}}$ in forms of $A_{\alpha}, B_{\beta}$,
the conclusion of the relative coefficients are same as above via similar arguments.
Without lose of generality, we can assume $g(m,n)\equiv0$.
From the above discussion we know, the coefficients defined is \eqref{ISM-eq-0} satisfy:
\begin{eqnarray}\label{ISM-eq-67}
  &&f(m,n)=\widetilde{f}(m,n)=0,\nonumber\\
  &&h(m,n)=\widetilde{h}(m,n)=0,\\
  &&g(m,n)=0,\ \widetilde{g}(m,n)=\widetilde{g}(m,n).\nonumber
\end{eqnarray}

{\bf Proof of theorem \ref{ISM-theorem main}:}
If $M=M_{\overline{0}}\oplus M_{\overline{1}}$ is the irreducible intermediate series module of $\mathfrak{L}$.
From the irreducibility and \eqref{ISM-eq-67}, we know $M_{\overline{1}}=0$.
This means  $I_{k}, G_{k}, H_{k}$ act trivially $M$, $\forall k\in\mathbb{Z}$.
This complete the proof of theorem \ref{ISM-theorem main}.
\qed

\section{Unitary representation of ${\mathfrak{L}}$}

Denote $\mathfrak{h}=\rm{Span}_{\mathbb{C}}\{L_{0}\}$ the Cartan subalgebra of $\mathfrak{L}$.

\subsection{Conjugate linear anti-involution of $\mathfrak{L}$}
Denote $\mathbf{A}=\rm{Span}_{\mathbb{C}}\{A_{k} \,|\,k\in\mathbb{Z}\}$,
in which $\mathbf{A}\in\{\mathbf{L}, \mathbf{I}, \mathbf{G}, \mathbf{H}\}$,
$A\in\{L,I, G, H\}$. We have the following lemma.

\begin{lemm}\label{UR-Lemma-Maximal ideal}
 $\mathbf{I}\oplus\mathbf{G}\oplus\mathbf{H}$\ is the unique maximal ideal of $\mathfrak{L}$.
\end{lemm}

\ni{\it Proof.} From the Lie super bracket of $\mathfrak{L}$, by direct calculation we know:
\begin{itemize}
  \item[(1)]$\mathbf{I}\oplus\mathbf{G}\oplus\mathbf{H}$ is an ideal of $\mathfrak{L}$.
  \item[(2)]None of $\mathfrak{L}$\'s ideal contain an nonzero element of $\mathbf{L}$.
\end{itemize}
This complete the proof of this lemma.

\begin{defi}\label{UR-Definition-anti-involution}
  $(1)$\ A map\ $\theta$: $\mathfrak{L}\rightarrow\mathfrak{L}$ is called
  the conjugate linear anti-involution of $\mathfrak{L}$\
  if $\theta$ satisfies the following conditions:
  \begin{enumerate}
    \item[C1:]\ \ $\theta(x+y)=\theta(x)+\theta(y)$,
    \item[C2:]\ \ $\theta(\alpha x)=\overline{\alpha}x$,
    \item[C3:]\ \ $\theta([x,y])=[\theta(y),\theta(x)]$,
    \item[C4:]\ \ ${\theta}^{2}=\rm{id}$.
  \end{enumerate}
  in which $x, y\in\mathfrak{L}, \alpha\in\mathbb{C}$,
  $\rm{id}$\ is the identity map of $\mathfrak{L}$.

  $(2)$ An $\mathfrak{L}$--module $M$ is called the unitary $\mathfrak{L}$--module,
  If there exists a positive defined  Hermitian form $\langle, \rangle$ on $M$
  such that $\forall u, v\in M, x\in\mathfrak{L}$:
  \[
    \langle xu,v\rangle=\langle u,\theta(x)v\rangle.
  \]
\end{defi}

The already known results about the conjugate linear anti-involution of the Virasoro algebra is listed in the following:
(one can see \cite{CP-CM1998}):
\begin{theo}\label{UR-Theorem-anti-involution of Vir}
  The conjugate linear anti-involution of the Virasoro algebra is one of the following.
  \begin{enumerate}
    \item[(1)] ${\theta}^{+}_{\alpha}(L_{m})={\alpha}^{m}L_{-m}$, in which $\alpha\in{\mathbb{R}}^{*}$.
    \item[(2)] ${\theta}^{-}_{\alpha}(L_{m})=-{\alpha}^{m}L_{m}$, in which $\alpha\in S^{1}$, the set of all length $1$ complex nunber .
  \end{enumerate}
\end{theo}

\begin{lemm}\label{UR-Lemma-theta-preserve CSA-MI}
  If $\theta$ is an arbitrary conjugate linear anti-involution of  $\mathfrak{L}$, we have:
  \begin{enumerate}
    \item[(1)] $\theta(\mathbf{I}+\mathbf{G}+\mathbf{H})=\mathbf{I}+\mathbf{G}+\mathbf{H}$.
    \item[(2)] $\theta(\mathfrak{h})=\mathfrak{h}$.
  \end{enumerate}
\end{lemm}

\ni{\it Proof.}  $(1)$  $\forall x\in\mathfrak{L}, y\in\mathbf{I}+\mathbf{G}+\mathbf{H}$,
From equation
\[
 [x,\theta(y)]=\theta[y,\theta(x)]
\]
and the fact $\mathbf{I}+\mathbf{G}+\mathbf{H}$ is an ideal of $\mathfrak{L}$, we have:
$\theta(\mathbf{I}+\mathbf{G}+\mathbf{H})$ is an ideal of $\mathfrak{L}$.
once more, we get:
\begin{equation}\label{UR-eq-1}
\theta(\mathbf{I}+\mathbf{G}+\mathbf{H})\subset \mathbf{I}+\mathbf{G}+\mathbf{H}.
\end{equation}
(If not,
From the ideal $\theta(\mathbf{I}+\mathbf{G}+\mathbf{H})$ we get another maximal ideal of
$\mathbf{I}+\mathbf{G}+\mathbf{H}$,
contradicts with the fact that $\mathbf{I}+\mathbf{G}+\mathbf{H}$ is the unique
maximal ideal of $\mathfrak{L}$).
\begin{equation}\label{UR-eq-2}
  {\theta}^{2}(\mathbf{I}+\mathbf{G}+\mathbf{H})\subset \theta(\mathbf{I}+\mathbf{G}+\mathbf{H}).
\end{equation}
Since ${\theta}^{2}=\rm{id}$, together with \eqref{UR-eq-1} and \eqref{UR-eq-2}, we have:
\[
 \theta(\mathbf{I}+\mathbf{G}+\mathbf{H})=\mathbf{I}+\mathbf{G}+\mathbf{H}.
\]

$(2)$ From $C3$ of the definition \ref{UR-Definition-anti-involution}, we have:
\[
 [\theta(L_{0}),\theta(L_{m})]=\theta([L_{m},L_{0}])=m\theta(L_{m}).
\]
\[
 [\theta(L_{0}),\theta(I_{m})]=\theta([I_{m},L_{0}])=m\theta(I_{m}).
\]
\[
 [\theta(L_{0}),\theta(G_{m})]=\theta([G_{m},L_{0}])=m\theta(G_{m}).
\]
\[
 [\theta(L_{0}),\theta(H_{m})]=\theta([H_{m},L_{0}])=m\theta(H_{m}).
\]
This implies $\theta(L_{0})$ acts diagonally on $\mathfrak{L}$.
We get $\theta(L_{0})\in\mathfrak{h}$.
Since $\theta$ is non-degenerate, together with $\rm{dim}(\mathfrak{h})=1$, we know:
$\theta(\mathfrak{h})=\mathfrak{h}$.
\qed

Suppose $\theta$ is an arbitrary conjugate linear anti-involution of $\mathfrak{L}$ , assume:
\begin{eqnarray}\label{UR-eq-formal sum of theta}
  \theta(L_{m})\!\!\!&=&\!\!\!\sum_{k}a_{m,k}^{L}L_{k}+\sum_{k}b_{m,k}^{L}I_{k}+\sum_{k}c_{m,k}^{L}G_{k}+\sum_{k}d_{m,k}^{L}H_{k},\nonumber\\
  \theta(I_{m})\!\!\!&=&\!\!\!\sum_{k}b_{m,k}^{I}I_{k}+\sum_{k}c_{m,k}^{I}G_{k}+\sum_{k}d_{m,k}^{I}H_{k},\nonumber\\
  \theta(G_{m})\!\!\!&=&\!\!\!\sum_{k}b_{m,k}^{G}I_{k}+\sum_{k}c_{m,k}^{G}G_{k}+\sum_{k}d_{m,k}^{G}H_{k},\\
  \theta(H_{m})\!\!\!&=&\!\!\!\sum_{k}b_{m,k}^{H}I_{k}+\sum_{k}c_{m,k}^{H}G_{k}+\sum_{k}d_{m,k}^{H}H_{k}.\nonumber
\end{eqnarray}
Only finitely many summands on the right side of the above four equations are non zero,
all the coefficients are lie in $\mathbb{C}$.

Since $\mathbf{I}+\mathbf{G}+\mathbf{H}$ is the maximal ideal of $\mathfrak{L}$,
and
\[
 {\mathfrak{L}}/{(\mathbf{I}+\mathbf{G}+\mathbf{H})}\cong Vir
\]
Here the notation Vir means the centerless  Virasoro algebra.
From theorem \ref{UR-Theorem-anti-involution of Vir} we know  $\theta(L_{n})$ must be one of the following two forms:
\begin{lemm}\label{UR-Lemma-formal sum of thetaL}
  \begin{enumerate}
    \item[(a)] ${\theta}^{+}_{\alpha}(L_{m})={\alpha}^{m}L_{-m}+\sum_{k}b_{m,k}^{L}I_{k}+\sum_{k}c_{m,k}^{L}G_{k}+\sum_{k}d_{m,k}^{L}H_{k}$,
    in which  $\alpha\in{\mathbb{R}}^{*}$.
    \item[(b)] ${\theta}^{-}_{\alpha}(L_{m})=-{\alpha}^{m}L_{m}+\sum_{k}b_{m,k}^{L}I_{k}+\sum_{k}c_{m,k}^{L}G_{k}+\sum_{k}d_{m,k}^{L}H_{k}$,
    in which  $\alpha\in S^{1}$, the set of all complex numbers of length $1$.
  \end{enumerate}
\end{lemm}

We determinate $\theta$ case by case.

{\bf Case $1$}\ \  $\theta(L_{n})$ has the form $(a)$ in \ref{UR-Lemma-formal sum of thetaL}.

From $(2)$ of Lemma \ref{UR-Lemma-theta-preserve CSA-MI} we know $\theta(L_{0})=L_{0}$.
Since
\begin{equation}\label{UR-eq-3}
 [\theta(L_{n}),\theta(L_{0})]=-n\theta(L_{n}),
\end{equation}
together with \eqref{UR-eq-formal sum of theta}, we have:
\begin{eqnarray}\label{UR-eq-4}
  &&n{\alpha}^{n}+n\sum_{k}b_{n,k}^{L}I_{k}+n\sum_{k}c_{n,k}^{L}G_{k}+n\sum_{k}d_{n,k}^{L}H_{k}\nonumber\\
  &&=n{\alpha}^{n}-\sum_{k}kb_{n,k}^{L}I_{k}-\sum_{k}kc_{n,k}^{L}G_{k}-\sum_{k}kd_{n,k}^{L}H_{k}
\end{eqnarray}
Compare the coefficients of the summands respectively in equation \eqref{UR-eq-4}, we have:
\[
 b_{n,k}^{L}=c_{n,k}^{L}=d_{n,k}^{L}=0,\ \forall\, k\neq-n.
\]
This implies:
\begin{equation}\label{UR-eq-5}
  \theta(L_{n})={\alpha}^{n}L_{-n}+b_{n,-n}^{L}I_{-n}+c_{n,-n}^{L}G_{-n}+d_{n,-n}^{L}H_{-n}.
\end{equation}
Similar arguments as \eqref{UR-eq-3}, for $\theta(A_{m})$ we have:
\begin{equation}\label{UR-eq-6}
  b_{m,k}^{A}=c_{m,k}^{A}=d_{m,k}^{A}=0,\ \forall\, k\neq-m.
\end{equation}
In which $A\in\{I, G, M\}$.
Since
\begin{equation}\label{UR-eq-7}
  (m-n)\theta(L_{m+n})=\theta([L_{m},L_{n}])=[\theta(L_{n}),\theta(L_{m})],
\end{equation}
together with \eqref{UR-eq-formal sum of theta},  replace \eqref{UR-eq-6} in the above equation,
compare the coefficients of each summands after the concrete calculation, we get:
\begin{equation}\label{UR-eq-8}
  (m-n)b_{m+n,-m-n}^{L}={\alpha}^{n}(m-n)b_{m,-m}^{L}+{\alpha}^{m}(m-n)b_{n,-n}^{L}+c_{m,-m}^{L}c_{n,-n}^{L}.
\end{equation}
\begin{equation}\label{UR-eq-9}
  (m-n)c_{m+n,-m-n}^{L}={\alpha}^{n}(m-\frac{n}{2})c_{m,-m}^{L}+{\alpha}^{m}(\frac{m}{2}-n)c_{n,-n}^{L}.
\end{equation}
\begin{eqnarray}\label{UR-eq-10}
  (m-n)d_{m+n,-m-n}^{L}\!\!\!&=&\!\!\!{\alpha}^{n}(m-\frac{n}{2})d_{m,-m}^{L}+(2m-n)b_{n,-n}^{L}c_{m,-m}^{L}\nonumber\\
  &&+(m-2n)b_{m,-m}^{L}c_{n,-n}^{L}+(\frac{m}{2}-n){\alpha}^{m}d_{n,-n}^{L}.
\end{eqnarray}
In \eqref{UR-eq-9}, take $m=n\neq0$, we know $c_{k,-k}^{L}=0$, $k\neq0$.
Again in \eqref{UR-eq-9}, take $m=-n=1$,  we get $c_{0,0}^{L}=0$.
Take together, we have:
\begin{equation}\label{UR-eq-11}
  c_{k,-k}^{L}=0,\ \forall\, k\in\mathbb{Z}.
\end{equation}
Now, \eqref{UR-eq-8} becomes:
\begin{equation}\label{UR-eq-12}
  (m-n)b_{m+n,-m-n}^{L}=(m-n)({\alpha}^{n}b_{m,-m}^{L}+{\alpha}^{m}b_{n,-n}^{L}).
\end{equation}
By induction on \eqref{UR-eq-12}:
\begin{equation}\label{UR-eq-13}
  b_{k,-k}^{L}=k{\alpha}^{k-1}b_{1,-1}^{L}.
\end{equation}
\ \eqref{UR-eq-10} becomes:
\begin{equation}\label{UR-eq-14}
  (m-n)d_{m+n,-m-n}^{L}={\alpha}^{n}(m-\frac{n}{2})d_{m,-m}^{L}+{\alpha}^{m}(\frac{m}{2}-n)d_{n,-n}^{L}.
\end{equation}
By induction on \eqref{UR-eq-14}:
\begin{equation}\label{UR-eq-15}
  d_{k,-k}^{L}=k{\alpha}^{k-1}d_{1,-1}^{L}.
\end{equation}
From  \eqref{UR-eq-11} \eqref{UR-eq-13} and \eqref{UR-eq-15}, we get:
\begin{equation}\label{UR-eq-16}
  \theta(L_{k})={\alpha}^{k}L_{-k}+k{\alpha}^{k-1}b_{1,-1}^{L}I_{-k}+k{\alpha}^{k-1}d_{1,-1}^{L}H_{-k}.
\end{equation}

From the super Lie bracket in $\mathfrak{L}$,
together with  (1) in lemma \ref{UR-Lemma-theta-preserve CSA-MI} and \eqref{UR-eq-6},
we can assume:
\begin{eqnarray}\label{UR-eq-17}
  \theta(I_{m})\!\!\!&=&\!\!\!b_{m,-m}^{I}I_{-m}+d_{m,-m}^{I}H_{-m},\nonumber\\
  \theta(G_{m})\!\!\!&=&\!\!\!b_{m,-m}^{G}I_{-m}+c_{m,-m}^{G}G_{-m}+d_{m,-m}^{G}H_{-m},\\
  \theta(H_{m})\!\!\!&=&\!\!\!d_{m,-m}^{H}H_{-m}.\nonumber
\end{eqnarray}

Plug \eqref{UR-eq-17} into
\[
 \theta([I_{n},L_{m}])=[\theta(L_{m}),\theta(I_{n})],
\]
compare the coefficients of each summands respectively, we have:
\begin{equation}\label{UR-eq-18}
  (n-m)b_{m+n,-m-n}^{I}=(n-m){\alpha}^{m}b_{n,-n}^{I},
\end{equation}
\begin{equation}\label{UR-eq-19}
  (n-m)d_{m+n,-m-n}^{I}=(n-\frac{m}{2}){\alpha}^{m}d_{n,-n}^{I}.
\end{equation}
By induction on \eqref{UR-eq-18}, we get:
\begin{equation}\label{UR-eq-20}
  b_{k,-k}^{I}={\alpha}^{k}b_{0,0}^{I},\ \forall\, k\in\mathbb{Z}.
\end{equation}
In \eqref{UR-eq-19}, take $m=n\neq0$:
\[
  d_{k,-k}^{I}=0,\ k\neq0.
\]
Again in \eqref{UR-eq-19}, take $m=-n=1$, we have:
\[
 d_{0,0}^{I}=0,
\]
this implies:
\begin{equation}\label{UR-eq-21}
  d_{k,-k}^{I}=0,\ \forall\, k\in\mathbb{Z}.
\end{equation}
Associative with \eqref{UR-eq-20} and \eqref{UR-eq-21}, we get:
\begin{equation}\label{UR-eq-22}
  \theta(I_{k})={\alpha}^{k}b_{0,0}^{I}I_{-k}.
\end{equation}

Plug \eqref{UR-eq-17} into
\[
 \theta([G_{n},L_{m}])=[\theta(L_{m}),\theta(G_{n})],
\]
compare the coefficients of each summands respectively, we have:
\begin{equation}\label{UR-eq-23}
  (n-\frac{m}{2})b_{m+n,-m-n}^{G}=(n-m){\alpha}^{m}b_{n,-n}^{G},
\end{equation}
\begin{equation}\label{UR-eq-24}
  (n-\frac{m}{2})c_{m+n,-m-n}^{G}=(n-\frac{m}{2}){\alpha}^{m}c_{n,-n}^{G},
\end{equation}
\begin{equation}\label{UR-eq-25}
  (n-\frac{m}{2})d_{m+n,-m-n}^{G}=(n-\frac{m}{2}){\alpha}^{m}d_{n,-n}^{G}+(2n-m)m{\alpha}^{m-1}b_{1,-1}^{L}c_{n,-n}^{G}.
\end{equation}
 \eqref{UR-eq-23}\ implies:
\begin{equation}\label{UR-eq-26}
  b_{k,-k}^{G}=0,\ \forall\, k\in\mathbb{Z}.
\end{equation}
By induction on \eqref{UR-eq-24}, we get:
\begin{equation}\label{UR-eq-27}
  c_{k,-k}^{G}={\alpha}^{k}c_{0,0}^{G}.
\end{equation}
By induction on \eqref{UR-eq-25}:
\begin{equation}\label{UR-eq-28}
  d_{k,-k}^{G}={\alpha}^{k}d_{0,0}^{G}+2k{\alpha}^{k-1}b_{1,-1}^{G}c_{0,0}^{G}.
\end{equation}
Associate \eqref{UR-eq-26}\ \eqref{UR-eq-27}\ \eqref{UR-eq-28}, we get:
\begin{equation}\label{UR-eq-29}
  \theta(G_{K})={\alpha}^{k}c_{0,0}^{G}G_{-k}+({\alpha}^{k}d_{0,0}^{G}+2k{\alpha}^{k-1}b_{1,-1}^{G}c_{0,0}^{G})H_{-k}.
\end{equation}

Plug \eqref{UR-eq-17} into
\[
 \theta([H_{n},L_{m}])=[\theta(L_{m}),\theta(H_{n})],
\]
compare the coefficients of each summands respectively, we have:
\begin{equation}\label{UR-eq-30}
  (n-\frac{m}{2})d_{m+n,-m-n}^{H}=(n-\frac{m}{2}){\alpha}^{m}d_{n,-n}^{H}.
\end{equation}
from this we know:
\begin{equation}\label{UR-eq-31}
  d_{k,-k}^{H}={\alpha}^{m}d_{0,0}^{H}.
\end{equation}
this implies:
\begin{equation}\label{UR-eq-32}
  \theta(H_{k})={\alpha}^{k}d_{0,0}^{H}H_{-k}.
\end{equation}

Plug \eqref{UR-eq-17} into
\[
 \theta([G_{m},G_{n}])=[\theta(G_{n}),\theta(G_{m})],
\]
compare the coefficients of each summands respectively, we have:
\begin{equation}\label{UR-eq-33}
  b_{0,0}^{I}={(c_{0,0}^{G})}^{2}.
\end{equation}

Plug \eqref{UR-eq-17} into
\[
 \theta([G_{n},I_{m}])=[\theta(I_{m}),\theta(G_{n})],
\]
compare the coefficients of each summands respectively, together together with \eqref{UR-eq-33}, we have:
\begin{equation}\label{UR-eq-34}
  d_{0,0}^{H}=b_{0,0}^{I}c_{0,0}^{G}={(c_{0,0}^{G})}^{3}.
\end{equation}

Associate \eqref{UR-eq-16},\ \eqref{UR-eq-22},\ \eqref{UR-eq-29},\
\eqref{UR-eq-32},\ \eqref{UR-eq-33},\ \eqref{UR-eq-24}, we have:
\begin{equation}\label{UR-eq-35}
  \begin{cases}
    \theta(L_{k})={\alpha}^{k}L_{-k}+k{\alpha}^{k-1}b_{1,-1}^{L}I_{-k}+k{\alpha}^{k-1}d_{1,-1}^{L}H_{-k},\\\vspace*{8pt}

    \theta(I_{k})={\alpha}^{k}{(c_{0,0}^{G})}^{2}I_{-k},\\\vspace*{8pt}

    \theta(G_{k})={\alpha}^{k}c_{0,0}^{G}G_{-k}+({\alpha}^{k}d_{0,0}^{G}+2k{\alpha}^{k-1}b_{1,-1}^{L}c_{0,0}^{G})H_{-k},\\ \vspace*{8pt}

    \theta(H_{k})={\alpha}^{k}{(c_{0,0}^{G})}^{3}H_{-k}.
  \end{cases}
\end{equation}
From \eqref{UR-eq-35} and C4 of the definition \ref{UR-Definition-anti-involution}, we get:
\begin{equation}\label{UR-eq-36}
  \begin{cases}
    \alpha\in{\mathbb{R}}^{*},\\\vspace*{8pt}

    c_{0,0}^{G}={\rm{e}}^{\rm{i}{\omega}_{1}},\ \omega_{1}\in[0,2\pi)\\\vspace*{8pt}

    b_{1,-1}^{L}=|b_{1,-1}^{L}|{\rm{e}}^{\rm{i}({\omega}_{1}+\delta\pi)},\\\vspace*{8pt}

    d_{1,-1}^{L}=|d_{1,-1}^{L}|{\rm{e}}^{\rm{i}(\frac{3}{2}{\omega}_{1}+\delta\pi)},\\\vspace*{8pt}

    d_{0,0}^{G}=|d_{0,0}^{G}|{\rm{e}}^{\rm{i}(2{\omega}_{1}+\delta\pi)},\\\vspace*{8pt}
    \delta\in\{0,1\}.
  \end{cases}
\end{equation}

{\bf Case\ $2$}\ \ $\theta(L_{n})$ has the form $(b)$ in \ref{UR-Lemma-formal sum of thetaL}.

Similar arguments as in case $1$, we have the following results:

\begin{equation}\label{UR-eq-37}
  \begin{cases}
    \theta(L_{k})=-{\alpha}^{k}L_{k}+k{\alpha}^{k-1}b_{1,1}^{L}+k{\alpha}^{k-1}d_{1,1}^{L}H_{k},\\\vspace*{8pt}

    \theta(I_{k})={\alpha}^{k}{c_{0,0}^{G}}^{2}I_{k},\\\vspace*{8pt}

    \theta(G_{k})={\alpha}^{k}c_{0,0}^{G}G_{k}+({\alpha}^{k}d_{0,0}^{G}-2k{\alpha}^{k-1}b_{1,1}^{L}c_{0,0}^{G})H_{k},\\\vspace*{8pt}

    \theta(H_{k})=-{\alpha}^{k}{(c_{0,0}^{G})}^{3}H_{k}.
  \end{cases}
\end{equation}
In which
\begin{equation}\label{UR-eq-38}
  \begin{cases}
    \alpha={\rm{e}}^{\rm{i}\omega_{2}},\ \omega_{2}\in[0,2\pi)\\\vspace*{8pt}

    c_{0,0}^{G}={\rm{e}}^{\rm{i}\omega_{3}},\ \omega_{3}\in[0,2\pi)\\\vspace*{8pt}

    b_{1,1}^{L}=|b_{1,1}^{L}|{\rm{e}}^{\rm{i}(\omega_{2}+\omega_{3}+\delta\pi)},\\\vspace*{8pt}

    d_{0,0}^{G}=|d_{0,0}^{G}|{\rm}^{\rm{i}(2\omega_{3}+\delta\pi)},\\\vspace*{8pt}

    d_{1,1}^{L}=|d_{1,1}^{L}|{\rm{e}}^{\rm{i}(\omega_{2}+\frac{3}{2}\omega_{3}+(\delta+1)\pi)},\\\vspace*{8pt}

    \delta\in\{0,1\}.
  \end{cases}
\end{equation}
\qed

From the discussion above, we have the following theorem:
\begin{theo}\label{UR-Theorem-theta}
  The conjugate linear anti-involution of $\mathfrak{L}$ must lies in one of the following forms:
  \begin{enumerate}
    \item[(1)] Given by \eqref{UR-eq-35} and \eqref{UR-eq-36}.
    \item[(2)] Given by \eqref{UR-eq-37} and \eqref{UR-eq-38}.
  \end{enumerate}
\end{theo}

\subsection{Irreducible intermediate series unitary module of $\mathfrak{L}$}
The unitary representation of the Virasoro algebra, owns the following already known results: (To be concrete, one can see \cite{CP-CM1998}).
\begin{theo}\label{UR-Theorem-URofVir}
 The irreducible unitary module of the Virasoro algebra is either the highest or lowest weight module,
 or isomorphism to the intermediate series module $A_{a,b}$.
 In which $a\in\mathbb{R}, b\in\frac{1}{2}+\rm{i}\mathbb{R}$.
\end{theo}

From the discussion of the irreducible intermediate series module of $\mathfrak{L}$ discussed in the previous subsection:
the irreducible intermediate series module $M$ of $\mathfrak{L}$ must also be the irreducible intermediate series module of the Virasoro algebra,
in other word, $\forall k\in\mathbb{Z}$, $I_{k}, G_{k}, H_{k}$ act trivially on $M$.
From the discussion above, we have the following theorem:

\begin{theo}\label{UR-Theorem-ISMURofL}
  The irreducible intermediate series unitary module of $\mathfrak{L}$
  is isomorphism to the irreducible intermediate series module of form
  $A_{a,b,0,0,0}$, in which $a\in\mathbb{R}, b\in\frac{1}{2}+\rm{i}\mathbb{R}$.
\end{theo}

\end{document}